\documentclass[12pt,a4paper]{amsart}

\usepackage{amscd}
\usepackage{graphicx}
\usepackage{amssymb,amsfonts,latexsym}
\usepackage{fixltx2e}
\usepackage{mathscinet}

\input xy
\xyoption {all}

\newtheorem{thm}{Theorem}[section]

\newtheorem{cor}[thm]{Corollary}

\newtheorem{lem}[thm]{Lemma}

\newtheorem{prop}[thm]{Proposition}

\newtheorem{conj}[thm]{Conjecture}

\theoremstyle{definition}

\newtheorem{defn}[thm]{Definition}

\newtheorem{ques}[thm]{Question}

\numberwithin{equation}{section}

\theoremstyle{remark}
\newtheorem{rem}[thm]{Remark}

\newcommand\oline[1] {{\overline{#1}}}

\def\divides{{\,|\,}}

\newcommand{\mM}{\mathcal{M}}

\newcommand{\fp}{\mathcal{\frak{p}}}
\newcommand{\ra}{\rightarrow}
\newcommand{\mQ}{\mathbb{Q}}
\newcommand{\mZ}{\mathbb{Z}}
\long\def\thesis#1\endt{#1}

\DeclareMathOperator\Gal{Gal}
\DeclareMathOperator\PSL{PSL}

\DeclareMathOperator\Hom{Hom}

\DeclareMathOperator\HLG{H}

\DeclareMathOperator\lcm{lcm}

\DeclareMathOperator\res{res}
\DeclareMathOperator\SL{SL}

\DeclareMathOperator\im{Im}
\DeclareMathOperator\coker{coker}

\begin{document}

\title%[] %% short title
{Tamely ramified subfields of division algebras}%

\def\UMich{Department of Mathematics, 530 Church St., University of Michigan, Ann Arbor 48109. }

\author{ Danny Neftin}
\address{\UMich}
\email{neftin@umich.edu}%

\begin{abstract}
For any number field $K$, it is unknown which finite groups appear as Galois groups of  extensions $L/K$ such that $L$ is a maximal subfield of a  division algebra with center $K$ (a $K$-division algebra). For $K=\mQ$, the answer is described by the long standing $\mQ$-admissibility conjecture.

We extend a theorem of Neukirch on embedding problems with local constraints in order to determine for every number field $K$, what finite solvable groups $G$ appear as Galois groups of tame maximal subfields of  $K$-division algebras, generalizing Liedahl's theorem for metacyclic $G$ and Sonn's solution of the $\mQ$-admissibility conjecture for solvable groups. %Under the further assumption that $G$ is metacyclic this was proven by Liedahl  and for $K=\mQ$ this proven

\end{abstract}

\maketitle

\section{Introduction}

A  division algebra $D$ which is finite dimensional over its center $K$ ({\it a $K$-division algebra}), is called a {\it $G$-crossed product} if there exists a Galois extension $L/K$ with Galois group $G$ (a {\it $G$-extension}) such that $L$ is a maximal subfield of $D$. Crossed products  are fundamental in the study of division algebras and are accompanied by a structure which explicitly describes them (see \cite[Chp. 14-19]{P}). % \cite{ABGV}).
A group $G$ is called {\it $K$-admissible} if there exists a $G$-crossed product $K$-division algebra;  a field extension $L/K$ is called {\it adequate} if $L$ is a maximal  subfield of a $K$-division algebra \footnote{In fact by \cite{Sch}, $L/K$ is adequate if and only if $L$ is a subfield of a $K$-division algebra. Thus, the maximality requirement can be omitted.}.

It is known by the Brauer-Hasse-Noether theorem that over a number field $K$,  all $K$-division algebras are crossed products with respect to a cyclic group. However, it is unknown for which groups $G$ there exists a $G$-crossed product $K$-division algebra, i.e. what groups are $K$-admissible?
%Crossed products are important!!!

Over $\mQ$, Schacher observed (\cite{Sch}) that the Sylow subgroups $P$ of a $\mQ$-admissible group are {\it metacyclic}, that is $P$ has a cyclic normal subgroup $C\lhd P$ such that $P/C$ is also cyclic.
The converse of this observation is known as the $\mQ$-admissibility conjecture: %It is conjectured that the converse also holds:
\begin{conj}\label{Q-ad.conj} Every group with metacyclic Sylow subgroups is $\mQ$-admissible. \end{conj}

This conjecture was studied extensively (e.g. \cite{FS},\cite{Fei},\cite{Fei2},\cite{FV},\cite{GS},\cite{Sch}) and  proven by Sonn for solvable groups in a series of papers (\cite{Son2}, \cite{Son1} and \cite{Son3}).
%\begin{thm}{\rm (Sonn)}\label{Q-solv.thm} Every solvable group with metacyclic Sylow subgroups  is $\mQ$-admissible.  \end{thm}

Recently, analogs of this conjecture were proved by Harbater, Hartmann and Krashen  over function fields of curves over complete discretely valued fields with algebraically closed residue fields (\cite{HHK}, cf. \cite{HH12}), by Paran and the author  over two dimensional complete local domains with algebraically closed residue fields (\cite{NP}), and by Surendranath and Suresh   over function fields of curves over complete discretely valued fields which contain enough roots of unity (\cite{SuSu}). However, the situation over number fields is far from being understood.

Schacher's observation extends to number fields under an additional assumption of tameness as follows. %was extended by Liedahl (see \cite[Theorem 28]{Lid2} or Remark \ref{Lid.rem}).  %A field extension $L/K$ is called {\it adequate} if $L$ is a maximal subfield of a finite dimensional $K$-division algebra.
Let $\mu_n$ denote the set of $n$-th roots of unity and $\sigma_{t,n}$ be the automorphism of $\mQ(\mu_n)$ for which $\sigma_{t,n}(\zeta)=\zeta^t$ for all $\zeta\in\mu_n$.
Using a similar argument to Liedahl's \cite[Theorem 28]{Lid2}, we observe that if $G$ appears as a Galois group of a  tamely ramified adequate extension  of a number field $K$ then its Sylow subgroups are metacyclic, and furthermore for every $l\divides |G|$,
%Over a general number field $K$,  Schacher's observation was extended by Liedahl (see \cite{Lid2}) as follows.
%If $G$ is $K$-admissible  and $p$ does not decompose in $K$ (there is only one prime of $K$ that divides $p$) then for every $p\divides |G|$
the $l$-Sylow subgroups $G(l)$ of $G$ admit a presentation: %are metacyclic and admit a presentation:
\begin{equation}\label{metacyclic.pre}
G(l)\cong \mathcal{M}(m,n,i,t):=\langle x,y| x^{m}=y^{i}, y^{n}=1, x^{-1}yx=y^{t}\rangle
\end{equation}
such that $\sigma_{t,n}$ fixes $K\cap \mQ(\mu_{n})$ (see ``only if part" of Theorem \ref{main.thm}).

 This observation suggests the following natural generalization of Conjecture \ref{Q-ad.conj}:
\begin{ques}\label{tame.prob} Let $K$ be a number field and $G$ a group whose $l$-Sylow subgroups admit a presentation $\mM(m,n,i,t)$ %as in {\rm (\ref{metacyclic.pre})}
such that $\sigma_{t,n}$ fixes $K\cap \mQ(\mu_{n})$, for every $l\divides |G|$. Is $G$ necessarily $K$-admissible? Furthermore, is there necessarily a tamely ramified adequate $G$-extension of $K$?
\end{ques}
The first part of this question is known to have an affirmative answer for  metacyclic $G$ (\cite[Theorem 27]{Lid2}) and for some small order groups: $A_5$ (\cite{GS}), the central extension $\SL_2(5)$ of $A_5$ (\cite{Fei5}), $A_6,A_7$ (\cite{SS}), the double covers of $A_6$ and $A_7$ (\cite{Fei4}),  $\PSL_2(7)$ (\cite{Al}) and  $\PSL_2(11)$ (\cite{Fei3}). %and for several small order non-solvable groups for $A_5$ and make a list.
%\begin{thm}{\rm (Liedahl)}\label{Lid.thm} Let $K$ be a number field and $G$ a metacyclic group such that for every $p||G|$ the $p$-Sylow subgroups of $G$ admit a presentation as in {\rm (\ref{metacyclic.pre})} such that $\sigma_{t,n}\in \Gal(\mQ(\mu_{n})/(\mQ(\mu_{n})\cap K))$. Then $G$ is $K$-admissible.
%\end{thm}%In fact, a similar argument (see Section \ref{prelim.sec}) shows Liedahl's observation can be used to prove that
%The adequate extensions which were constructed in theorems \ref{Q-solv.thm} and \ref{Lid.thm} are all tamely ramified.

 In this paper we give a positive answer to Question \ref{tame.prob} for solvable groups, generalizing Liedahl's \cite[Theorem 27]{Lid2} and Sonn's \cite[Theorem 1]{Son3}: % by studying  tamely ramified adequate extensions. % and use them to solve Problem \ref{tame.prob} affirmiatively for solvable groups: %the following generalization of theorems \ref{Q-solv.thm} and \ref{Lid.thm}:
% In fact, we prove:
\begin{thm}\label{main.thm}
Let $K$ be a number field and $G$ a solvable group. Then there exists a tamely ramified adequate $G$-extension $L/K$ if and only if for every $l||G|$, the $l$-Sylow subgroups of $G$ admit a presentation $\mM(m,n,i,t)$ such that $\sigma_{t,n}$ fixes $K\cap \mQ(\mu_n)$. %\mQ(\mu_{n})\cap K)).$$
\end{thm}

%The proof of Sonn's theorem over $\mQ$ is based on two theorems of Neukirch, \cite[Main Theorem]{Neu} and \cite[Korollar 6.4]{Neu2}.
We note that since the proof of Sonn's theorem (\cite{Son3}) over $\mQ$ is based on Neukirch's \cite[Main Theorem]{Neu} which makes an assumption on the absence of roots of unity in $K$,  Sonn's proof does not apply over arbitrary number fields.

A key ingredient in our proof is an extension of \cite[Korollar 6.4]{Neu2}. Neukirch's Korollar 6.4 is a highly useful tool that  under the assumption of at least one of six conditions on a finite set $S$ of primes of the base field, allows to change solutions of embedding problems to satisfy any prescribed local conditions at $S$ (generalizing the Grunwald-Wang theorem). We extend Korollar 6.4 by showing that under the assumption of at least one of four of these six conditions on $S$, it is possible to change a solution to satisfy  prescribed conditions at $S$ leaving the solution unchanged at any given finite set of primes $T$.

We use this extension to strengthen Sonn's proof of \cite[Theorem 1]{Son3} in order to obtain tamely ramified adequate $G$-extensions of $\mQ$ with prescribed local behavior at given finite sets of primes. This gives us a strong control over the ramification  of $G$-crossed product $\mQ$-division algebras, allowing us  to lift these to division algebras over a given number field and by that prove Theorem \ref{main.thm}.

%We also provide  a reduction of Question \ref{tame.prob} for non-solvable groups, using \cite{Son2}, to a realization problem with prescribed local condition of a list of six families of groups.
 % generalize Sonn's theorem (\cite{Son3}) by obtaini

%We take the approach of realizing $G$ over $\mQ$ by adequate extensions and lifting them to $K$. However,

% Note that Sonn's proof of Theorem \ref{Q-solv.thm} does not extend directly to arbitrary number fields since it relies on a theorem of Neukirch \cite[Main Theorem]{Neu} that assumes $K$ does not have any non-trivial odd order roots of unity.
%Credit to Jack by saying the proof of the strengthening is based on his proof.
%Our method
%%In order to prove Theorem \ref{main.thm} our method will be
%is to prove a refinement of Theorem \ref{Q-solv.thm} that yields many $\mQ$-adequate $G$-extensions $L/\mQ$ and prove that for some of these extensions $LK/K$ is a $K$-adequate $G$-extension. The proof of the refinement is an adaptation of Sonn's proof of Theorem \ref{Q-solv.thm}.

%\subsection*{Acknowledgments}

%--------------------------------- Embedding problems  --------------------------------

This work is partially based on the author's Ph.D. thesis (\cite{Nef}). %and Theorem \ref{main.thm} is an improvement of \cite[Theorem 0.1.6]{Nef}.
I would like to thank my thesis advisor Jack Sonn for investing time and effort into teaching and guiding me,  and for helpful comments on this manuscript.

\section{Embedding problems and local Galois groups}\label{embedding.sec}

\subsection{Embedding problems}
%The main ingredients of our proof arise from
The theory of embedding problems  is central in the study of the inverse Galois problem and is a key ingredient in our proof of Theorem \ref{main.thm}.  We shall  describe a setup for these problems, recall Neukirch's \cite[Korollar 6.4]{Neu2} and extend it.

\subsubsection{Setup} Embedding problems are a strong generalization of the inverse Galois problem which ask whether a Galois extension can be embedded into a larger Galois extension with a given Galois group. The precise setup is as follows.

 A  {\it (finite) embedding problem} over a number field $K$ consists of a finite Galois extension $L/K$ and an epimorphism of finite groups $\pi:E\ra G:=\Gal(L/K)$.  For our purposes it suffices to consider embedding problems with abelian kernel $A:=\ker(\pi)$.

Let $G_K$ denote the absolute Galois group of $K$.
Two homomorphisms $\psi_1,\psi_2:G_K\ra E$ are called equivalent if there is an $a\in A$  such that $a^{-1}\psi_1(g)a=\psi_2(g)$ for all $g\in G_K$. A {\it solution} for $\pi$ is an equivalence class of homomorphisms $\psi:G_K\ra E$ (not necessarily surjective)
%that makes diagram (\ref{section1-embedding problem}) commutative, i.e
for which $\pi\circ\psi$ is the restriction map $\res_L:G_K\ra G$.
For a surjective solution $\psi$, the  fixed field %\footnote{The field $M$ is sometimes called a Galois cover of $L/K$}
$M=\oline K^{\ker(\psi)}$  contains $L$ and has Galois group $\Gal(M/K)\cong E$.

The epimorphism $\pi$ defines an action of $G$ on $A$ and hence induces a $G_K$-module structure on $A$  via  $\res_L$.
For every crossed homomorphism $\chi\in\HLG^1(G_K,A)$ and solution $\psi:G_K\ra E$,  the map $\psi'=\chi\cdot\psi$ given by $\psi'(\sigma)=\chi(\sigma)\psi(\sigma)$ for all $\sigma\in G_K$, is also a solution (see \cite[Chp. IX, \S 4]{NSW}).
In fact, for every two solutions $\psi,\psi'$ of $\pi$, there is a unique $\chi\in\HLG^1(G_K,A)$ such that $\psi'=\chi\cdot\psi$.
We think of $\chi$ as the element that ``changes" $\psi$ to $\psi'$.

\subsubsection{Embedding problems with prescribed local conditions}
By a prime $\fp$ of $K$ we mean a finite or infinite prime. Fix an algebraic closure $\oline{K}$ of $K$, an algebraic closure $\oline{K}_\fp$ of the completion $K_\fp$, and an inclusion of $\oline{K}$ into $\oline{K}_\fp$ for every prime  $\fp$ of $K$. In particular,
the embedding problem $\pi$ induces a local embedding problem $\pi_\fp:\pi^{-1}(G_\fp)\ra G_\fp$ where $G_\fp=\Gal(L_\fp/K_\fp),L_\fp:=LK_\fp$. %, for every prime $\fp$.
%For a prime $\fp$ of $K$,
Moreover, the restriction $\psi_\fp$ of a solution $\psi:G_K\ra E$ to the subgroup $G_{K_\fp}$ is a solution for $\pi_\fp$.

Let $S$ be a finite set of primes of $K$ and for every $\fp\in S$  fix (prescribe) a solution $\psi^{(\fp)}$  to  $\pi_\fp$, assuming such (local) solutions exist. Similarly to the Grunwald-Wang theorem,  one is interested in solutions $\psi$ of $\pi$ such that $\psi_\fp=\psi^{(\fp)}$ for all $\fp\in S$.

Assume $\pi$ has a solution $\phi$.
Then for every $\fp\in S$ there is  $\chi^{(\fp)}\in \HLG^1(G_{K_\fp},A)$ such that $\psi^{(\fp)}=\chi^{(\fp)}\cdot \phi_\fp$. If the element $(\chi^{(\fp)})_{\fp\in S}$ has a source $\chi$ under the restriction map:
\begin{equation*}\label{changes.equ} \rho_S:\HLG^1(G_{K},A)\ra \prod_{\fp\in S}\HLG^1(G_{K_\fp},A)\end{equation*}
then $\psi:=\chi\cdot\phi$ is a solution for $\pi$ which restricts to $\psi^{(\fp)}=\chi^{(\fp)}\cdot\phi_\fp$ at all  $\fp\in S$. Thus, if the map $\rho_S$ is surjective, every solution for  $\pi$  can be ``changed" to a solution with  prescribed local conditions at $S$.

\subsubsection{Neukirch's Korollar} \cite[Korollar 6.4]{Neu2} is  a highly useful criteria for the map $\rho_S$ to be surjective. Let $A$ be a $G_K$-module and $n=\exp(A)$.  Let $A'=\Hom(A,\mu_n)$ be the dual $G_K$-module  %with the action $f^\sigma(a)=f(a^{\sigma^{-1}})^\sigma$ for all $a\in A,\sigma\in G_K$ and $f\in A'$. Let $$G_K'=\{\sigma\in G_K| f^\sigma=f\mbox{ for all } f\in A'\},$$
and  $K(A')$  the fixed field of the centralizer of $A'$ in $G_K$. %$G_K'$ and
Let $G':=\Gal(K(A')/K)$ and for a prime $\fp$ of $K$, let $G_\fp':=\Gal(K(A')_\fp/K_\fp)$. %, where $w$ is a prime divisor of $v$ in $K(A')$. The group $G_v'$ is well defined since for a different prime divisor $w'$ of $v$ one has $$\Gal(K(A')_{w'}/K_v)\cong \Gal(K(A')_w/K_v).$$%\footnote{I don't think this is what Neukirch ment}
%Since the action of  $G_K'$ on $A'$ is trivial, the action of $G_K$ on $A'$ induces an action of $G'$ on $A'$. % and hence an action of $G_\fp'$ on $A'$. Let
Denote
%$$ \Gamma(G_\fp',A'):=\ker\left(\HLG^1(G_\fp',A')\ra \prod_{g\in G_\fp'}\HLG^1(\langle g\rangle,A')\right).$$
$ \Gamma(G,A):=\ker\left(\HLG^1(G,A)\ra \prod_{g\in G}\HLG^1(\langle g\rangle,A)\right).$

\begin{thm}\emph{(Neukirch \cite[Korollar 6.4]{Neu2})}\label{7cond.thm} %Let $K$ be a number field and $S$ a finite set of primes of $K$. Let $A$ be a finite $G_K$-module.
Let $S$ be a finite set of primes of  $K$. Then the map
%\begin{equation}
%\HLG^1(G_K,A)\ra \prod_{v\in S} \HLG^1(G_{K_v},A)
%\end{equation}
$\rho_S$ is surjective in each of the following cases:
\begin{enumerate}
\item[\rm (a)] $\Gamma(G_\fp',A')=0$ for all $\fp\in S$,
\item[\rm (b)] for every $\fp\in S$, the group $G_\fp'$ is cyclic or a semidirect product of two cyclic groups of relatively prime orders,
\item[\rm (c)] $\HLG^1(G',A')=0$,
\item[\rm (d)] $|G'|=\lcm\{|G_\fp'||\fp\not\in S\}$,
\item[\rm (e)] $A$ is cyclic of odd order,
\item[\rm (f)] the action of $G_K$ on $A$ is trivial and $(K,\exp(A),S)$ does not fall into a special case. %Definition \ref{special.defn}).
\end{enumerate}
In {\rm (f)}, when $\exp(A)=2^tm$, $m$ odd,  one says that the triple $(K,\exp(A),S)$  falls into a special case if $K(\mu_{2^t})/K$ is noncyclic and $S$ contains all primes $\fp$ for which $K_\fp(\mu_{2^t})/K_\fp$ is noncyclic.
\end{thm}

%
%Let $K$ be a number field and $M/K$ a Galois extension with Galois group $G$. Let $\pi:E\ra G$ be an embedding problem with abelian kernel $A$ that has a (weak) solution $\psi$.  Let $S,W$ be finite sets of primes of $K$.
%
%Following Neukrich (ref) we let $A'=\Hom(A,\mu)$ be the dual $G_K$-module, $K(A')$ the fixed field of the centralizer of $A'$ in $G_K$, $G':=\Gal(K(A')/K)$, $G'_v:=\Gal(K(A')_v/K_v)$ and
%$$ \Gamma(G,A):=\ker ( \HLG^1(G, A)\ra \prod_{\sigma\in G} \HLG^1(\langle\sigma\rangle, A)). $$
%
%Under the assumption $\Gamma(G_v',A')=0$ for all $v\in S$, Korollar 6.4a) of (ref) states that the restriction map:
%$$ res_S:\HLG^1(G_K,A)\ra \prod_{v\in S} \HLG^1(G_{K_v},A) $$
%is surjective. We show that furthermore:

%Thus, under any of these conditions a global solution guarantees a solution with prescribed local conditions.
%Conditions (a)-(b) on the other hand  depend on $S$ and may hold for $S$ but fail for $S\cup T$ for some finite set $T$.

Thus, under each of these conditions one can change a solution to satisfy arbitrary prescribed local conditions at $S$. Furthermore, we show that under each of conditions (a), (b), (c) or (e) it is possible to change a solution to satisfy prescribed local conditions at $S$ leaving the solution unchanged at a given finite set of primes~$T$.
%We shall be interested in changing a solution for $\pi$ to have prescribed restrictions at $S$ without changing its restrictions to a given finite set of primes  $T$. We show that this is possible if conditions (a) and (b) hold for $S$ even though they may fail for $S\cup T$.
%However, we will show that conditions (a) or (b) insure more than the surjectivity of $\rho_S$. Namely, under these conditions one can change a solution for $\pi$  to satisfy arbitrary local conditions at $S$ leaving the restrictions at a finite set $T$ unchanged.
%In fact, if one of the conditions (a-c) or (e) holds then has: %one can also preserve the solution at a finite set of primes $S$:
\begin{prop}\label{neu.gen} Let $A$ be a finite $G_K$-module. %Let  $S$ and $T$ be two disjoint finite sets of primes of $K$.
Assume that conditions (a) or (b) hold for a finite set $S$. Then the subgroup
$$\prod_{\fp\in S}\HLG^1(G_{K_\fp},A)\times \prod_{\fp\in T}\{0\}$$ is in the image of $\rho_{S\cup T}$ for every finite set $T$ disjoint from $S$.
\end{prop}

\begin{proof} %Conditions (c) and (e) are independent of $S$ and therefore show that $\res_{S\cup T}$ is in fact surjective.
Since by \cite[Satz 6.2]{Neu2} condition (b) implies (a), it suffices to prove the assertion when (a) holds.
Assume that $\Gamma(G'_\fp,A')=0$ for all $\fp\in S$. Let $P$ be the set of all primes of $K$ and $\prod_{\fp\in P}'\HLG^1(G_{K_\fp},A)$ the restricted product over the subgroup $\prod_{\fp\in P}\HLG^1_{un}(G_{K_\fp},A)$.
Recall that the Poitou-Tate theorem gives a non-degenerate bilinear map
$$ \beta:\prod_{\fp\in P}\stackrel{'}{} \HLG^1(G_{K_\fp},A)\times \prod_{\fp\in P}\stackrel{'}{}\HLG^1(G_{K_\fp},A') \ra \mQ/\mZ $$
which is defined as the product of local bilinear maps $$\beta_\fp:\HLG^1(G_{K_\fp},A)\times \HLG^1(G_{K_\fp},A')\ra \mQ/\mZ$$ for every $\fp\in P$. % in which the images of $\HLG^1(G_{K},A)$ and $\HLG^1(G_{K},A')$ are orthogonal complements.

Following  \cite{Neu2}, for a finite set $U$ of primes of $K$ we let:
%$$\rho_{U}:\HLG^1(G_K,A)\ra \prod_{\fp\in U}\stackrel{'}{} \HLG^1(G_{K_\fp},A)$$
$$\rho_{U}':\HLG^1(G_K,A')\ra \prod_{\fp\not\in U}\stackrel{'}{} \HLG^1(G_{K_\fp},A')$$
be the restriction map, $\Delta=\coker(\rho_{S\cup T}) $ and $\nabla=\ker(\rho_{S\cup T}')/\ker(\rho'_\emptyset)$.
%For an element $x$ in (resp. subset $A$ of) a product over a set of primes $U$, we denote by $x_\fp$ (resp. $A_\fp$) the $\fp$-th component of the product (resp. the projection to the $\fp$-th component).
By \cite[Satz 4.4]{Neu2}, %the Poitu-Tate theorem
$\beta$ induces a non-degenerate bilinear form
$\beta_0:\Delta\times~\nabla\ra\mQ/\mZ$, which is given on $\chi:=(\chi_\fp)_{\fp\in S\cup T} \in \Delta$ and $\lambda\in \nabla$ by $\beta_0(\chi,\lambda):=\beta(\tilde\chi,\rho'_\emptyset(\lambda))$ where $\tilde\chi\in \prod_{\fp\in P}'\HLG^1(G_{K_\fp},A)$ is any  element
%Let $\tilde\chi\in \prod_{\fp\in P}\HLG^1(G_{K_\fp},A)$ be %an element (a lift) which is
%$0$ at all $\fp\not\in S\cup T$ and
whose $\fp$-th component is $\chi_\fp$ at all $\fp\in S\cup T$. % Then $\beta_0(\chi,\lambda)=\beta(\tilde\chi,\rho'_\emptyset(\lambda))$.

Let $\chi=\prod_{\fp\in S\cup T}\chi_\fp$ be an element of $\Delta$ such that $\chi_\fp=0$ for all $\fp\in T$. We claim that $\chi$ is orthogonal to $\nabla$ and therefore it is the zero element in $\Delta$, proving the proposition.

Letting $\tilde\chi= (\tilde\chi_\fp)_{\fp\in P}\in \prod_{\fp\in P}'\HLG^1(G_{K_\fp},A)$ where $\tilde\chi_\fp=\chi_\fp$ for $\fp\in S$ and  $\tilde\chi_\fp=0$ for $\fp\not\in S$,
we have $\beta_0(\chi,\nabla)=\beta(\tilde\chi,\rho_\emptyset'(\nabla))$.
Since $\tilde\chi_\fp=0$ for  $\fp\not\in S$, it suffices to show that $\beta_\fp(\chi_\fp,\rho'_\emptyset(\nabla)_\fp)=0$ for all $\fp\in S$, where $\rho'_\emptyset(\nabla)_\fp$ is the projection of $\rho'_\emptyset(\nabla)$ to the $\fp$-th factor. %$\HLG^1(G_{K_v},A)$.
But \cite[Satz 6.3]{Neu2} implies that the image of $\nabla$ under the restriction map
$$ \rho_{S\cup T,A'}:\HLG^1(G_K,A') \ra \prod_{\fp\in S\cup T} \HLG^1(G_{K_\fp},A') $$
lies in $\prod_{\fp\in S\cup T}\Gamma(G_\fp',A')$.
Since by assumption $\Gamma(G_\fp',A')=0$ for $\fp\in S$, we get  $\rho_\emptyset'(\nabla)_\fp=\rho_{S\cup T,A'}(\nabla)_\fp=0$ and hence $\beta_\fp(\chi_\fp,\pi_\fp\rho'_\emptyset(\nabla))=0$ for all $\fp\in S$, proving the claim.
%$\rho'_{S\cup T}$ is injective. %$\pi_\fp\rho'_\emptyset
% We therefore claim that $\pi_\fp\rho_\emptyset'(\nabla)=0$ for all $\fp\in S$.

%As in Satz 6.3, let $\oline\rho_{S\cup T}$ be the restriction map $\HLG^1(G',A')\ra \prod_{\fp\not\in S\cup T}'\HLG^1(G_\fp',A')$.
%The restriction maps give  a commutative diagram:
%\begin{equation} \xymatrix{
%                   \ker(\oline\rho_s') \ar[r]^{\pi_\fp\oline\rho'_\emptyset} \ar[d] &  \HLG^1(G_\fp',A') \ar[d]   \\
%                    \ker(\rho_s') \ar[r]^{\pi_\fp\rho'_\emptyset}   & \HLG^1(G_{K_\fp},A') \\
%                   } \end{equation}
%for every $\fp\in S$. But as shown in Satz 6.3, the left vertical map is an isomorphism and the image of $\pi_\fp\oline\rho_\emptyset'$ lies in  $\Gamma(G_\fp',A')$ which is assumed to be $0$. Therefore $\pi_\fp\rho_\emptyset'(\ker(\rho_s'))=0$, proving the claim.
%
%%For $v\in S$, let $\chi_v'\in \HLG^1(G_{K_v},A)$ be a lift of $\chi_v\in
\end{proof}

%As a consequence, we can change solutions to a prescribed solution at $S$ so that they remains the same at a finite set of primes $T$: %

From Proposition \ref{neu.gen} and the discussion above it we get:
\begin{cor}\label{Neu.cor} Let $\pi:E\ra \Gal(L/K)$ be an embedding problem with  solution $\phi$.
%Let $S$
%and $T$ be finite disjoint sets of primes and
Fix  solutions $\psi^{(\fp)}$ for $\pi_\fp$ at all primes $\fp$ in a  finite set $S$ and let $T$ be a finite set of primes disjoint from $S$. Assume that at least one of conditions {\rm (a),(b),(c)} or {\rm (e)} hold for $S$. %$\Gamma(G_\fp',A')=0$ for all $\fp\in S$.

Then there exists a solution $\psi$ such that $\psi_\fp=\psi^{(\fp)}$ for all $\fp\in S$ and $\psi_\fp=\phi_\fp$ for all $\fp\in T$.
\end{cor}
\begin{proof} Since conditions (c) and (e)  are independent of $S$, the image of $\rho_{S\cup T}$ contains  %and hence insure that $\rho_S$ is surjective for every $S$.
$\prod_{\fp\in S}\HLG^1(G_{K_\fp},A)\times \prod_{\fp\in T}\{0\}$ under these conditions as well. For $\fp\in S$, let $\chi^{(\fp)}\in \HLG^1(G_{K_\fp},A)$ be the element for which $\psi^{(\fp)}=\chi^{(\fp)}\cdot \psi_\fp$. By Proposition \ref{neu.gen}, the element $(\chi^{(\fp)})_{\fp\in S}\times  (0)_{\fp\in T}$ has a source $\chi\in \HLG^1(G_K,A)$ under the map $\rho_{S\cup T}$. Then the solution $\psi:=\chi\cdot\phi$ restricts to $\chi^{(\fp)}\cdot \phi_\fp=\psi^{(\fp)}$ at all $\fp\in S$ and to $0\cdot \phi_\fp=\phi_\fp$ at all $\fp\in T$.
\end{proof}
\begin{rem}\begin{enumerate} \item Proposition \ref{neu.gen} need not hold under conditions (d) or (f). For example, let $K$ be a quadratic extension of $\mQ$ in which $2$ splits and let $\fp_1,\fp_2$ be the primes above it. Let $S=\{\fp_1\}$, $T=\{\fp_2\}$ and let $A=\mZ/8$ be the trivial $G_K$-module. Then $A'\cong \mu_8$ as  $G_K$-modules and $K(A')=K(\mu_8)$.  Both conditions (d) and (f) hold for $S$ and hence $\rho_S$ is surjective.

However, since  $K(A')_\fp/K_\fp$ is cyclic for all $\fp\neq\fp_1,\fp_2$,  conditions (d) and (f) fail for $S\cup T$.
Furthermore, the Grunwald-Wang theorem shows that $$ \prod_{\fp\in S}\HLG^1(G_{K_{\fp}},A)\times \prod_{\fp\in T} \{0\}\not\subseteq \im\rho_{S\cup T}. $$
Indeed, letting  $\psi^{(\fp_2)}=0$ and $\psi^{(\fp_1)}\in\Hom(G_{K_{\fp_1}}, A)$ be such that the fixed field of $\ker(\psi^{(\fp_1)})$ is the unramified $\mZ/8$-extension of $K$, \cite[Chp. X, Theorem 5]{AT} shows that  $(\psi^{(\fp_1)},\psi^{(\fp_2)})\not\in \im(\rho_{S\cup T}).$
%On the other hand, both conditions (d) and (f) hold for $S$ and hence $\rho_S$ is surjective.
\item Let $A$ be a trivial $G_K$-module of exponent $2^tm'$ where $m'$ is odd. If $K(\mu_{2^t})/K$ is cyclic then condition (f) holds for all finite sets $S$. %and in particular if $\mu_{2^t}\subseteq K$.
\end{enumerate}
\end{rem}
%Discuss Grunwald Wang
\subsection{Tame Galois groups of local fields}\label{tame.sec}
We shall make use of a few well known facts about Galois groups of tame local extensions,  all of which can be found in \cite{Wei} and \cite[Chp. VII, \S 5]{NSW}.
 Let $L/K$  a tamely ramified $G$-extension %, i.e. a finite Galois extension with Galois group $G$.
of $p$-adic fields, $I$ its  inertia group, $n:=|I|$, and $q$ the cardinality of the residue field of $K$. %and  $q$ the cardinality of the residue field of $K$.

The subfield $L^I$ contains $\mu_n$ and $L/L^I$ is a (cyclic) Kummer extension.
The  Galois group of $L^I/K$ is generated by the Frobenius automorphism $\sigma_L$ which acts on $\mu_n$ by raising each element to the power  $q$. In particular the restriction of $\sigma_L$ to $\mQ(\mu_n)$ is $\sigma_{q,n}$ and fixes $\mQ(\mu_n)\cap K$.  The action of $\sigma$ on $I$ via  conjugation in $G$ is equivalent to its action on $\mu_n$. Thus, %  by raising each element to the power of $q$.  In particular, $\Gal(L/K)$ is isomorphic to:
\begin{equation*}%\label{presentation.equ}
%\begin{array}{c}
G= \mM(m,n,i,t)=  \langle x,y| x^{m}=y^{i}, y^{n}=1, x^{-1}yx=y^{t}\rangle,  \\ %\end{equation}
% \mbox{ and } \sigma_{q,n}\in \Gal(\mQ(\mu_{n})/(\mQ(\mu_{n})\cap K))
%\end{array}
\end{equation*}
where $t\equiv q$ (mod $n$), $I=\langle y\rangle$ and $x=\sigma$ (mod $I$). In particular, one has the following observation which is the basis to \cite[Theorem 28]{Lid2}:
\begin{lem}\label{Liedahl.lem} Let $p,l$ be two distinct rational primes and $K$  a  $p$-adic field.  Then every group $G$ that appears as a Galois group over $K$ has a metacyclic $l$-Sylow subgroup $\mM(m,n,i,t)$ such that $\sigma_{t,n}$ fixes $K\cap \mQ(\mu_n)$. %\in\Gal(\mQ(\mu_n)/\mQ(\mu_n)\cap K)$.
\end{lem}
\begin{proof} Let $L/K$ be a $G$-extension, $M$ the fixed subfield of an $l$-Sylow subgroup $H$ of $G$ and $t$ the cardinality of the residue field of $M$. Then  $H\cong\mM(m,n,i,t)$ % where $t\equiv q$ (mod $n$).
and %In particular, $\sigma_{q,n}$ is the Frobenius automorphism of $\Gal(M(\mu_n)/M)$ (REF). Thus,
$\sigma_{t,n}$ fixes $M\cap \mQ(\mu_n)$.  In particular $\sigma_{t,n}$ fixes $K\cap \mQ(\mu_n)$.
\end{proof}
%Conversely, if $q\equiv t$ (mod $n$), then $\mu_n$ is contained in the (unique) unramified $m$-extension $M/K$ and

Consider the converse problem of realizing $\mM(m,n,i,t)$ over $K$ and assume $t\equiv q$ (mod $n$) so that $\mM(m,n,i,t)=\mM(m,n,i,q)$. The Galois group $G^{tr}_K$ of the maximal tamely ramified extension of $K$ is profinitely generated by two automorphisms $\sigma$ and $\tau$ and one relation $\sigma^{-1}\tau\sigma=\tau^q$, where $\tau$ is of order prime to $q$ and $\sigma$ is the Frobenius automorphism.  Letting $M$ be the (unique) unramified degree $m$ extension of $K$, %the restriction map $G_K\ra \Gal(M/K)$  splits through the restriction $G_K^{tr}$ and
 $\sigma$ restricts to the Frobenius automorphism of $M/K$. Thus, an embedding problem $\pi:\mM(m,n,i,t)\ra \Gal(M/K)$ with kernel $\langle y\rangle$ has a surjective solution whose corresponding field is a tamely ramified $\mM(m,n,i,t)$-extension of~$K$.
%Consider the converse direction of realizing a group $G=\mM(m,n,i,t)$ over $K$.
%\begin{enumerate}
% \item[(T1)]
% The extension $L/L^I$ is a (cyclic) totally ramified Kummer extension.  In particular, $\mu_n\subseteq L^I$ where $n=|I|$.
%\item[(T2)] The  Galois group of $L^I/K$ is generated by the Frobenious automorphism $\sigma$ which acts on $I$ by raising each element to the power of $q$. %, where $q$ is cardinality of the residue field of $K$.
 %and $m=[L^I:K]$.   (the number of roots of unity).
%Unramified extension followed by a totally ramified extension. The unramified extension must have the roots of unity and it acts on the inertia group by raising each element to the power of the norm.
% In particular, $\Gal(L/K)$ is isomorphic to:
%\begin{equation}\label{presentation.equ}
%\begin{array}{c} \mM(m,n,i,q):=  \langle x,y| x^{m}=y^{i}, y^{n}=1, x^{-1}yx=y^{q}\rangle,  \\ %\end{equation}
% \mbox{ and } \sigma_{q,n}\in \Gal(\mQ(\mu_{n})/(\mQ(\mu_{n})\cap K))
%\end{array}
%\end{equation}
%since the restriction of $\sigma$ to $\mQ(\mu_n)$ is $\sigma_{q,n}$.
%\item[(T3)]
%The embedding problem $\pi:\mM(m,n,i,q)\ra\Gal(M/K)$ where $M/K$ is the (unique) unramiifed $m$-extension and $\ker\pi=\langle y\rangle$, admits a surjective solution.
%\end{enumerate}

%%%----------------------------- Section on admissibility  ---------------------------------------------

\section{Galois groups of tamely ramified adequate extensions}\label{prelim.sec} %Let us introduce a stronger notion of admissibility

%Let $K$ be a number field.
\subsection{Proof of Theorem \ref{main.thm}}

%In the proof of Theorem \ref{main.thm}
We consider a refined notion of adequacy. For a number field $K$ and a finite set $S$ of primes of  $K$, we say that $L/K$ is {\it $S$-adequate} if $L$ is a maximal subfield of a $K$-division algebra that is  unramified outside $S$.  Let $D(L/K,\frak{P})$ denote the decomposition group of a prime $\frak{P}$ of $L$. The same proof as of Schacher's criterion (\cite[Proposition 2.6]{Sch}) gives the following criterion for $S$-adequacy:
%The basic criterion for adequacy over number fields~is:
%In \cite{Sch}, Schacher gave a necessary and sufficient criterion on a
%$G$-extension $L/K$   over a number field $K$ was given in \cite{Sch}. Namely, Schacher showed
%For a prime $\frak{P}$
\begin{prop}\label{Sch.thm}%{\rm (Schacher \cite{Sch})}
Let $L/K$ be a $G$-extension of number fields and $S$ a finite set of primes of $K$. Then $L/K$ is $S$-adequate if and only if
%A finite group $G$ is $K$-admissible if and only if
%there exists a Galois extension $L/K$ with Galois group $\Gal(L/K)\cong G$ such that %that satisfies:
%\emph{1}.  $\Gal(L/K) \cong G$,
%\emph{2}.
%that $L$ is $K$-adequate if and only if
for every rational prime $l\divides |G|$, there are two primes $\fp_1 ,\fp_2\in S$ for which  $D(L/K,\frak{P}_i)$ %$\Gal(L_{\frak{P}_i}/K_{\fp_i})$
contains an $l$-Sylow subgroup of $G$, where $\frak{P}_i$ is a prime of $L$ which divides $\fp_i$, $i=1,2$. % for $w_i|v_i$ and $i=1,2$.
\end{prop}
%Here $D(\frak{P}_i/\fp_i)$ denotes the decomposition group of $\frak{P}_i$ in $L/K$ and
Note that the condition of containing an $l$-Sylow subgroup of $G$ is independent of the choice of prime $\frak{P}_i$ dividing $\fp_i$.

A key ingredient in our proof of Theorem \ref{main.thm} is  the following generalization of Sonn's theorem (\cite[Theorem 1]{Son3}) which asserts the existence of $S$-adequate $G$-extensions for prescribed sets~$S$.
Since $\mM(m,n,i,t)$ is realizable over $\mQ_p$ when $p\equiv t$ (mod $n$) (see Section \ref{tame.sec}), we consider the following sets $S$:
%natural candidates $S$ for  $S$-adequacy of $G$-extensions: %extension with given Galois group: %over~$\mQ$:
 %supports for $G$:
\begin{defn}
%Let $G(l)$ be an $l$-Sylow subgroup of $G$.
We call a set $S$  of distinct odd rational primes $p_i^{(l)}, i=1,2,  l||G|$ which are prime to  $|G|$, a \emph{tame
supporting set  for $G$} if for every $l||G|$,  the $l$-Sylow subgroups of $G$ admit a presentation $\mathcal{M}(m,n,i,t)$  such that $p_1^{(l)}, p_2^{(l)}\equiv t$ (mod $n$). % and the primes of $S$ are prime to $|G|$.
\end{defn}

%In the following theorem we show that in fact, for every such set $S$, there exists an $S$-adequate $G$-extension of $\mQ$, generalizing Sonn's theorem (\cite{Son3}).

%We show that in fact Sonn's theorem (\cite{Son3}) holds much more generally by proving that for every such set $S$, there exits an $S$-adequate $G$-extension of $K$.

%We shall say that a set of primes $T$ {\it splits completely} in a field $L$ if every prime in $T$  splits completely in $L$.

%We shall call a group {\it Sylow metacyclic} if it has metacyclic Sylow subgroups.
%We shall adopt the proof of Theorem \ref{solvable sylow meta cyclic} to prove the following refinement:
%In fact, we show that if $G$ is solvable for every such set there exists an $S$-adequate $G$-extension of $\mQ$:
%We say that a set $T$ splits completely in a field $L$ if every prime in $T$ splits completely in $L$.
\begin{thm}\label{Q.thm}
Let $G$ be a solvable group with metacyclic Sylow subgroups. Let $S$ be a tame supporting set  for $G$  and $T$ a finite set of rational primes which is disjoint from $S$. %so that $S\cap T=\emptyset$. %(3) $v_i(2)\equiv 1$ (mod $u$) for every odd $u\in S\cup\{v_i(p)| p||G|,p\not=2,3\}$.
Then there exists an $S$-adequate   $G$-extension $L/\mQ$ in which the primes of $T$ split completely.
\end{thm}
The proof of this theorem is based on Corollary \ref{Neu.cor} and on Sonn's proof of \cite{Son3}, and is given in Section~\ref{refined.sec}.
We now use this theorem  to prove Theorem~\ref{main.thm}:
%\subsubsection{Splitting of primes}
%The proof of theorems \ref{main.thm} and \ref{Q.thm} also
%Theorem \ref{Q.thm} gives adequate extension in which  finitely many given primes split completely.

\begin{proof}[Proof of Theorem \ref{main.thm}] ``Only if part": Let $l$ be a prime that divides $|G|$. By Proposition \ref{Sch.thm}, there is a prime $\fp$  of $K$ such that $D(L/K,\frak{P}), \frak{P}\divides \fp,$ contains an $l$-Sylow subgroup $P$ of $G$. If $\fp\divides l$, $\fp$ is unramified in $L$ and hence $P$ is cyclic. Otherwise, Lemma \ref{Liedahl.lem} implies that $P$ has a presentation $\mM(m,n,i,t)$ such that $\sigma_{t,n}$ fixes  $K_\fp\cap \mQ(\mu_n)$ and hence $K\cap \mQ(\mu_n)$. In both cases, $P$ has the required presentation.

%\begin{proof}
``If part":  Let $l\divides |G|$ be a rational prime and let $P^{(l)}$ be the set of primes of $K$ that are  unramified over $\mQ$ with residue degree $1$, and whose restriction $p$ to $\mQ$ satisfies $p\equiv t_l$ (mod $n_l$).

%Frobenius automorphism in $K(\mu_n)/K$ is $\tilde{\sigma}_{t_l,n_l}$.

We first claim that $P^{(l)}$ is infinite. Let $M$ be the $\mQ$-normal closure of $K$. % and %$l\divides |G|$ a rational prime.
%$M(\mu_{n_l})$  whose Frobenius automorphism in $M(\mu_{n_l})/\mQ$ is
Since $\sigma_{t_l,n_l}$ %\in
fixes $K\cap \mQ(\mu_{n_l}),$ $\sigma_{t_l,n_l}$ extends to an automorphism %$\tilde{\sigma}_{t_l,n_l}$
of $K(\mu_{n_l})$ that fixes $K$ and hence lifts to an automorphism $\tau_l\in\Gal(M(\mu_{n_l})/K)$. % a lift of $\tilde{\sigma}_{t_l,n_l}$.
By  Chebotarev's density theorem there are infinitely many primes $\frak{P}$ of $M(\mu_{n_l})$ that are unramified over $\mQ$, and whose Frobenius automorphism in $M(\mu_{n_l})/\mQ$ is $\tau_l$. In particular, the restriction $p$ of such $\frak{P}$ to $\mQ$ has Frobenius $\sigma_{t_l,n_l}$ in $\mQ(\mu_{n_l})/\mQ$, and hence $p\equiv t_l$ (mod $n_l$). Since $\tau_l$ fixes $K$, the restriction of each such $\frak{P}$ to $K$ has residue degree $1$ over $\mQ$ and hence is in $P^{(l)}$, proving the claim.   %is ${\sigma}_{t_l,n_l}$.

Since $P^{(l)}$ is infinite, %. Thus, for every $l\divides |G|$,
we can choose two primes
$\fp_1^{(l)},\fp_2^{(l)}\in P^{(l)}$  such that the restrictions   $p_i^{(l)}$ of $\fp_i^{(l)}$ to $\mQ$, $i=1,2, l\divides |G|$, are all distinct rational primes which are prime to $|G|$.
%Therefore there is an automorphism $\tau_l\in \Gal(M(\mu_{n_l})/K)$ that  restricts to $\sigma_{t_l,n_l}$.
%Since $\sigma_{t_l,n_l}$ is the restriction of $\tau_l$ to $\mQ(\mu_{n_l})$, it is the Frobenius automorphism of $p_i^{(l)}$ in $\mQ(\mu_{n_l})$ and hence  $p_i^{(l)}\equiv t_l\mbox{ (mod }n_l\mbox{), for }i=1,2, l\divides |G|.$
Thus, the set $S:=\{p_i^{(l)}|i=1,2,l||G|\}$ is a tame supporting set for $G$ and
by Theorem \ref{Q.thm} there exists an $S$-adequate $G$-extension $L/\mQ$ in which all of the primes $l$ dividing $|G|$ split completely. %: whose primes split completely in $K$. %Let $\fp_i^{(l)}$ be a prime of $K$ which divides $p_i^{(l)}$ for $i=1,2, l\divides |G|$.

We claim that $N:=LK$ is an adequate extension of $K$. This proves the theorem since  $L/\mQ$ and hence $N/K$ is tamely ramified.  Since $\fp_i^{(l)}$ has residue degree $1$ over $\mQ$, we have $K_{\fp_i^{(l)}}\cong \mQ_{p_i}$ and hence:
\begin{equation}\label{degrees.equ} [N_{\frak{P}_i^{(l)}}:K_{\fp_i^{(l)}}]=[L_{\frak{P}_i^{(l)}\cap L}:\mQ_{p_i^{(l)}}]
 \mbox{ for  }\frak{P}_i^{(l)}\divides \fp_i^{(l)}, i=1,2,l\divides |G|. \end{equation} Letting $l^{\alpha_l}$ be the maximal power of $l$ dividing $|G|$, (\ref{degrees.equ}) shows that $l^{\alpha_l}\divides [N:K]$ for every $l\divides |G|$ and hence that $\Gal(N/K)\cong G$. Furthermore, (\ref{degrees.equ}) shows that
%By Theorem \ref{Q.thm} there exists an $S$-adequate $G$-extension $L/\mQ$. Since $p_i^{(l)}$ splits completely in $K$, Lemma \ref{disjoint.lem} shows that $L$ is disjoint from $K$ and hence $N:=LK$ is a $G$-extension of $K$. Furthermore, %letting  $\fp_i^{(l)}$ be a prime of $K$ which divides $p_i^{(l)}$ %for $i=1,2, l\divides |G|$.
%we have $[N_{\frak{P}_i^{(l)}}:K_{\fp_i^{(l)}}]=[L_{\frak{P}_i^{(l)}\cap L}:\mQ_{p_i^{(l)}}]$ for $\fp_i^{(l)}\divides p_i^{(l)}$ and $\frak{P}_i^{(l)}\divides \fp_i^{(l)}$, and hence
%$\Gal(N_{\frak{P}_i^{(l)}}/K_{\fp_i^{(l)}})$
$D(\frak{P}_i^{(l)}, N/K)$ contains an $l$-Sylow subgroup of $G$ for $i=1,2, l\divides |G|$, showing that $N/K$ is adequate, as required. %Note that $L/\mQ$ and hence $N/K$ are tamely ramified.
\end{proof}

%\begin{rem} A similar idea was used by Liedahl (\cite{Lid2}) to show the if part of Th
%\end{rem}
\begin{rem}\label{Lid.rem} \begin{enumerate} \item  In \cite{Lid2}, Liedahl uses Lemma \ref{Liedahl.lem}, similarly to the ``only if part" of Theorem \ref{main.thm}, to  show that under the assumption that $l$ does not decompose in  $K$, the $l$-Sylow subgroups of a $K$-admissible group admit a presentation $\mM(m,n,i,t)$ such that $\sigma_{t,n}$ fixes $K\cap \mQ(\mu_n)$. He also uses the flexibility of \cite[Theorem 1]{Son1} to prove  that if $G$ itself has a presentation $\mM(m,n,i,t)$ such that $\sigma_{t,n}$ fixes $K\cap \mQ(\mu_n)$, then $G$ is $K$-admissible.  %one has: % if $l$ does not decompose in $K$ then:
% \begin{equation}\label{Liedahl.equ} \begin{array}{l} \mbox{If $G$ is $K$-admissible then its $l$-Sylow subgroups  admit a presentation } \\ \mM(m,n,i,t)\mbox{ such that } \sigma_{t,n}\mbox{ fixes } K\cap \mQ(\mu_n), %(\ref{metacyclic.pre}),}
 %\end{array}\end{equation}  %under the assumption that $l$ does not decompose in  $K$,
% and the converse of (\ref{Liedahl.equ}) when $G$ is itself a metacyclic $l$-group.
%\end{rem}
%\begin{rem}
\item Note that the proof of the ``only if part" of Theorem \ref{main.thm} applies more generally without the assumption that $G$ is  solvable.
\item The proof of Theorem \ref{Q.thm} gives furthermore that $l^{\alpha_l}\divides [L_{\frak{P}_i^{(l)}\cap L}:\mQ_{p_i^{(l)}}]$ for all $l\divides |G|,i=1,2$.
\end{enumerate}
\end{rem}

\subsection{Consequences}\label{conclusions.sec}
%The following are consequences of Theorem \ref{main.thm}:
For solvable groups $G$ we get the following characterization of $K$-admissibility under the assumption that every $l\divides |G|$ does not decompose~in~$K$:
\begin{cor}\label{section2 - equivalence theorem}
Let $K$ be a number field and $G$  a solvable group. Assume that every prime $l$ that divides $|G|$ does not decompose in  $K$. Then the following conditions are equivalent:
\begin{enumerate}
%\item $G$ is tamely $K$-admissible,
\item  There exists a tamely ramified $K$-adequate $G$-extension; % $L/K$ which is everywhere tamely ramified,
\item $G$ is $K$-admissible;
%\item $G$ is $K$-preadmissible,
\item  for every $l||G|$, the $l$-Sylow subgroups of $G$ admit a presentation  $\mathcal{M}(m,n,i,t)$
    such that $\sigma_{t,n}$ fixes $K\cap \mQ(\mu_n)$. %\in\Gal(\mQ(\mu_n)/ \mQ(\mu_n)\cap K)$.
    %that satisfies (\ref{section3-Liedahls rep in thm on tame ad}). %for which \begin{equation*} \sigma_{t_p,n_p}\in \Gal(\mQ(\mu_{n_p})/(\mQ(\mu_{n_p})\cap K)). \end{equation*}
%\end{enumerate}
%If moreover $G$ is of odd order the conditions above are equivalent to:
%\begin{enumerate}
\end{enumerate}
\end{cor}
\begin{proof}
The implication $(1)\Rightarrow (2)$ is immediate and $(2)\Rightarrow (3)$ follows from Remark \ref{Lid.rem} (\cite[Theorem 28]{Lid2}). %$(6) \Rightarrow (1) \Rightarrow (2) \Rightarrow (3) $ are clear from the definitions and Lemma  \ref{restriction of k admissibile metacyclic groups}. The implication %$(5)\Rightarrow (1)$ follows directly from Theorem \ref{theroem tame admissibility of solvable groups}. We are left to prove
The implication $(3)\Rightarrow (1)$ is the ``if part" of Theorem~\ref{main.thm}. %for groups $G$ of odd order. Let $W$ be the set of primes of $K$ that whose restriction to $\mQ$ divides $|G|$. By Theorem \ref{theroem tame admissibility of solvable groups} there is a $K$-adequate Galois $G$-extension $L/K$ in which all primes of   $W$ split completely. In particular, $L/K$ is tamely ramified.
\end{proof}

%Under the assumptions of Theorem \ref{main.thm} the group $G$ is in fact infinitely often
%$K$-admissible (cf. \cite{AS}).
Recall that a group $G$ is called {\it infinitely often $K$-admissible}  if there exist infinitely many adequate $G$-extensions $L_i/K$, $i\in\mathbb{N}$, such that  $L_{r+1}\cap (L_1\cdots L_r)=~K$ (cf. \cite{Al}).
\begin{cor}\label{infinitely.cor} Let $K$ be a number field and $G$ a solvable group such that for every $l||G|$, the $l$-Sylow subgroups admit a presentation $\mathcal{M}(m,n,i,t)$ such that $\sigma_{t,n}$ fixes $K\cap \mQ(\mu_n)$. % (\ref{section3-Liedahls rep in thm on tame ad}) holds.
Then $G$ is  infinitely often $K$-admissible.

Furthermore, there exists a $K$-division algebra $D$ that has infinitely many disjoint maximal subfields $L_i,i\in\mathbb{N},$ such that $\Gal(L_i/K)\cong G$.
\end{cor}

The following lemma is useful in the proof of Theorem \ref{Q.thm} and is used here to prove Corollary \ref{infinitely.cor}. %and  produce adequate extensions that are disjoint from a given number field. %property of
Given a Galois extension $N/\mQ$, define the following condition on a finite set of rational primes~$T$:

% \begin{array}{cl}
\begin{enumerate} \item[(A$_N$)]  The decomposition groups $D(N/\mQ,\fp), p\in T, \fp\divides p, $ generate $\Gal(N/\mQ)$.
\end{enumerate}
%\mbox{ where $p\in T$ and $\fp$ is a prime of $M$ dividing $p$, %} \\
%&
%\mbox{ generate }\Gal(M/\mQ).$$
%\end{array}

%We shall use the following lemma to produce adequate extensions that are disjoint from a given number field.

\begin{lem}\label{disjoint.lem}
Assume $T$ satisfies $(A_N)$. Then every finite extension $K/\mQ$ in which the primes of $T$ split completely is disjoint from $N$. %, i.e. $L\cap M=K$. %If $G$
\end{lem}

%\begin{lem}\label{disjoint.lem}  Let $G$ be a group, $K$ a number field and $M/K$ a finite extension. If $G$ satisfies Condition \textup{\ref{section2- Condition for disjointness}} over $K$ then there is a $K$-adequate $G$-extension $L/K$ for which $L\cap M=K$. In particular $G$ is strongly $K$-admissible.
%\end{lem}

\begin{proof}
%Since $M$ and $L$ are Galois over $\mQ$, so is $M\cap L$.
%Let $K/\mQ$ be a finite extension in which the primes of $T$ split completely and
Let $H:=\Gal(N/N\cap K)$ and assume on the contrary that $H\neq G$. By condition (A$_N$) there exists a prime $p\in T$ and $\fp\divides p$  such that $D:=D(N/\mQ,\fp)\not\subseteq H$. In particular, $[N_\fp:K_{\fp\cap K}]=|D\cap H|<|D|=[N_\fp:\mQ_p]$ and hence $[K_{\fp\cap K}:\mQ_p]>1$ contradicting the assumption that $p$ splits completely in $K$.
%the decomposition group of $p$ in $M\cap L$ is $D_pH/H$ and hence $p$ does not split completely in $M\cap L$, contradicting the assumption that $p$ splits completely in $L$.
\end{proof}
Note that by Chebotarev's density theorem for every cyclic subgroup $C\leq G$ there are infinitely many primes $\fp$ of $N$ for which $D(N/\mQ,\fp)=C$. Thus, we can always choose a finite set $T$ which satisfies $(A_N)$.

\begin{proof}[Proof of Corollary \ref{infinitely.cor}]
% ===================== Proof without same division algebra ========================================
%It suffices to show that given  $r$ disjoint  $G$-extensions $L_1,...,L_r$ of $K$ there exists an adequate $G$-extension  $L_{r+1}/K$  such that $L_{r+1}\cap (L_1\cdots L_r)=K$.

%Let $S$ and $T$ be  as in the proof of Theorem \ref{main.thm}.  Let $M$ be the $\mQ$-normal closure of $K$ and $N:=L_1\cdot\cdot\cdot L_rM$.
%Extend $T$ to a finite set $T_0$  which is disjoint from $S$ and for which condition (A$_N$) holds.
%By Theorem \ref{Q.thm}, there exists an $S$-adequate $G$-extension $L/\mQ$ in which the primes of $T_0$ split completely. Thus, as in the proof of Theorem \ref{main.thm}, $L_{r+1}:=LK$ is an adequate $G$-extension of $K$. By Lemma \ref{disjoint.lem}, $L \cap N=\mQ$ and hence $L_{r+1}\cap %(L_1\cdots L_r)=K$.

% ==== Proof with the same division algebra =================
Let $S$ and $l^{\alpha_l}$ be  as in the proof of Theorem \ref{main.thm}. Define $D:=D_0\otimes_\mQ K$ where $D_0$ is the $\mQ$-division algebra with  Hasse invariants $1/{l^{\alpha_l}}$ at $p_1^{(l)}$, $-1/{l^{\alpha_l}}$ at $p_2^{(l)}$ for  $l\divides |G|$, and $0$ at all other primes.
%Let $D_0$ be as in Remark \ref{invariant.rem}.

It suffices to show that given  $r$ disjoint  $G$-extensions $L_1,...,L_r$ of $K$ which are maximal subfields of $D$ there exists a maximal subfield $L_{r+1}$ of $D$ such that $\Gal(L_{r+1}/K)=G$ and $L_{r+1}\cap (L_1\cdots L_r)=K$.

   Let $M$ be the $\mQ$-normal closure of $K$ and $N:=L_1\cdot\cdot\cdot L_rM$.
Let $T$ be a finite set  which is disjoint from $S$, contains all primes $l\divides |G|$, and for which condition (A$_N$) holds.

As remarked in \ref{Lid.rem}.(3), Theorem \ref{Q.thm} gives a maximal subfield  $L$ of $D_0$ in which the primes of  $T$ split completely.
By Lemma \ref{disjoint.lem}, $L\cap N=\mQ$ and hence $L_{r+1}:=LK$ is a $G$-extension of $K$ and $L_{r+1}\cap (L_1\cdots L_r)=K$. Since in addition $L_{r+1}$ splits $D$, $L_{r+1}$ is a maximal subfield of $D$.
%==============================================================
\end{proof}

%============== Remarking the same division algebra
%\begin{rem} In fact the fields $L_r,r\in \mathbb{N}$, in Corollary \ref{infinitely.cor} can all be chosen to be maximal subfields of the same $K$-division algebra $D$. Indeed, letting $l^{\alpha_l}$ be  the largest $l$-power dividing $|G|$, we define $D:=D_0\otimes_\mQ K$ where $D_0$ is the $\mQ$-division algebra with  Hasse invariants $1/{p^{\alpha_l}}$ at $p_1^{(l)}$, $-1/{p^{\alpha_l}}$ at $p_2^{(l)}$ for  $l\divides |G|$, and $0$ at all other primes. Since in the proof of Corollary \ref{infinitely.cor}, the fields $L_r$ can all be chosen to be of the form $L_r=LK$ where $L$ is an $S$-adequate $G$-extension obtained using an application of Theorem \ref{Q.thm}, Remark \ref{Lid.rem} shows that such fields $L_r$, $r\in\mathbb{N}$, are all maximal subfields of~$D$.
%\end{rem}
%============================

%-------------------------------- Generalization Section ------------------------------------------------------

\section{Proof of Theorem \ref{Q.thm}}\label{refined.sec}

In this section we prove Theorem \ref{Q.thm}. In Sections \ref{2-group.sec} and \ref{2-3-groups.sec} we treat the cases of $2$-groups and  $\{2,3\}$-groups (groups of order $2^a3^b$), respectively. We first show how the theorem follows from the latter case.

As in Section \ref{embedding.sec}, we fix an embedding of an algebraic closure of $\mQ$ into an algebraic closure of each of its completions.
We shall say that a set of primes $T$ {\it splits completely} in $L$ if every prime in $T$ splits completely in $L$.
%then deduce the general case.
%We can now %use \cite{Neu} to
%deduce Theorem \ref{Q.thm}.
\begin{proof}[Proof of Theorem \ref{Q.thm}]
Let $n=|G|$. By  \cite[Lemma 1.4]{Son2}, %$G$ has a $\{2,3\}$-normal complement. In other words,
there is a normal subgroup $N\lhd G$ of order prime to $2$ and $3$ and a $\{2,3\}$-subgroup $H$ such that $G=NH$.

Extend $T$ to a finite set $T_0$ disjoint from $S$ which satisfies condition $(A_{\mQ(\mu_n)})$.
%By Section \ref{2-3-groups.sec}, %Proposition \ref{section2- ++splits completely proposition of 2,3}
By the case of $\{2,3\}$-groups (Section \ref{2-3-groups.sec}),
there exists a $\{p_1^{(2)},p_2^{(2)},p_1^{(3)},p_2^{(3)}\}$-adequate $H$-extension $K/\mQ$ in which $T_0\cup \{p_1^{(l)},p_2^{(l)}|\, l> 3 \}$ splits completely. Since by Lemma \ref{disjoint.lem},  $K\cap \mQ(\mu_n)=\mQ$,
%Since $(|N|,|A|)=1$,
and since the embedding problem $G\ra \Gal(K/\mQ)$ splits, we may apply    \cite{Neu}.   %\ref{section1- Neukirch main embedding Theorem}.
It follows that  $K/\mQ$ embeds into a $G$-extension $L/\mQ$ such that  $T$ splits completely in $L$ and $\Gal(L_{p_i^{(l)}}/\mQ_{p_i^{(l)}})$, $i=1,2$, is an $l$-Sylow subgroup of $N$ for all $l||N|$. In particular, $L/\mQ$ is an $S$-adequate $G$-extension   in which $T$ splits completely, as required.
\end{proof}

\subsection{$2$-groups}\label{2-group.sec}

%\begin{prop}\label{section2- lemma of metacyclic 2-group} Let $G$ be a metacyclic $2$-group, $T$ a finite set of  rational primes and $S=\{p_1(2),p_2(2)\}$ a tame supporting set for $G$  disjoint from $T$. %for which $S\cap T=\emptyset$.  %, for which $v_i(2)\not= 3$, $i=1,2$.
%Then there exists a  $G$-extension that is compatible with $S$ and in which  $T$ splits completely.
%\end{prop}
%We use Corollary \ref{Neu.cor} to generalize
%Our proof for $2$-groups is based on
%The following is  a strong generalization of the argument in
The $\mQ$-admissibility of metacyclic $2$-group was proved in \cite{Son1} using Theorem \ref{7cond.thm}. We use Corollary \ref{neu.gen} in order to prove Theorem \ref{Q.thm} for $2$-groups, generalizing \cite{Son1}:
%The proof of \cite[Theorem 1]{Son1} is based on Theorem \ref{7cond.thm} and is generalized here by using Corollary \ref{Neu.cor}:
\begin{proof}[Proof of Theorem \ref{Q.thm} for $2$-groups]

Let $G$ be a metacyclic $2$-group with presentation $G\cong \mM(m,n,i,t)$ such that $\sigma_{t,n}$ fixes $K\cap \mQ(\mu_n)$,
 and let $k$ be the order of $x$ in $G$. Let $S=\{p_1,p_2\}$ be a tame supporting set for $G$ such that $p_i\equiv t$ (mod $n$),  $i=1,2$.

Since $S$ consists of odd primes, the Grunwald-Wang theorem (see Theorem \ref{7cond.thm}.(f)) implies that there exists a $\mZ/k$-extension $\hat{K}/\mQ$ in which the primes of $S$ are inert and $T$ splits completely. We identify $\Gal(\hat K/\mQ)$ with $\langle x\rangle$ and let $K/\mQ$ be the unique $\mZ/m$-extension inside $\hat{K}$. The embedding problem $\pi:G\ra \Gal(K/\mQ)$  with kernel $A:=\langle y\rangle$ has a solution $\phi:G_\mQ\ra \langle x\rangle\subseteq G$ which is given by the restriction map to $\Gal(\hat{K}/\mQ)$.  %We would like to change $\phi$ to a solution with the desired properties.

%At the prime $p_i$, $\pi$ induces a
Let $\pi_i: G\ra \Gal(K_{p_i}/\mQ_{p_i})$ be the corresponding local embedding problem at $p_i$, $i=1,2$.  Since $p_i\equiv t$ (mod $n$), $\pi_i$ has a surjective solution $\psi^{(i)}:G_{\mQ_{p_i}}\ra G$ whose fixed field $L^{(i)}$ is totally ramified over $K_{p_i}$  and in particular $\mu_n\subseteq K_{p_i}$, for $i=1,2$ (see Section \ref{tame.sec}).

In order to  change $\phi$ to a solution with the desired properties,  we apply Corollary \ref{Neu.cor}. Let $A'=\Hom(A,\mu_n)$ be the dual $G_\mQ$-module,
$K':=\mQ(A')$, $G'=\Gal(K'/\mQ)$ and  $G'_{p_i}:=\Gal(K'_{p_i}/\mQ_{p_i})$, $i=1,2$. %In particular $K'\subseteq K(\mu_n)$.
Since every automorphism in $G_\mQ$ that fixes $A$ and $\mu_n$ also fixes $A'$, we have $K'\subseteq K(\mu_n)$ and hence
$ K'_{p_i}\subseteq %(K(\mu_n))_{p_i}=
K_{p_i}(\mu_n)=K_{p_i},  $  for $i=1,2$.
Thus, $G'_{p_i}$ is cyclic and  %For $p\in W$, one has $k'_v\leq (k(\mu_n))_v=k_v(\mu_n)=\mQ_v(\mu_n)$. As $n$ is a $2$-power and $2\not\in W$, $\Gal(\mQ_v(\mu_n)/\mQ_v)$ and hence $G'_v$ are cyclic for all $v\in W$.
 condition \ref{7cond.thm}.(b) holds. By Corollary \ref{Neu.cor}, there exists a solution $\psi:G_\mQ\ra G$ of $\pi$, whose restriction at $p_i$ is  $\psi^{(i)}$, $i=1,2$,  and  the restriction  remains the trivial solution at each $p\in T$. Since $\psi^{(1)},\psi^{(2)}$ are surjective, $\psi$ is also surjective. Thus, the fixed field $L$ of $\ker\psi$ is an $S$-adequate $G$-extension of $\mQ$  in which $T$ splits completely, as required.
%Thus, $G'_v$ is cyclic for every $v\in W'$ and we may apply Theorem \ref{7cond.thm}.(b).

\end{proof}
\begin{rem}
%\begin{enumerate}
%\item  In the definition of a tame supporting set we assumed the set consists of odd primes in order to avoid the special case in the above application of the Grunwald-Wang theorem. %Theorem \ref{Q.thm} does not necessarily hold if $2\in S$. For example,  $G=\mZ/8\times \mZ/8$, the Grunwald-Wang theorem shows that there is no extension  .
%\item
Note that we use Corollary \ref{neu.gen} since  Theorem \ref{7cond.thm} cannot be applied for the set $S\cup T$. In fact, the Grunwald-Wang theorem shows that the map $\rho_{S\cup T}$ need not be surjective if $2\in T$.
%\end{enumerate}
\end{rem}

% ---------------------- Solvable Sylow metacyclic -----------------------------------
\subsection{$\{2,3\}$-groups}\label{2-3-groups.sec}

%%================ 2-3-proofs========================
%\subsubsection{$2$-Sylows with normal complement}
Let $G$ be a $\{2,3\}$-group and $G(3)$ a $3$-Sylow subgroup of $G$. If $G(3)$ is normal in $G$ then Theorem \ref{Q.thm}
 essentially follows from the $2$-groups case by applying \cite{Neu}:

\begin{proof}[Proof of Theorem \ref{Q.thm} for $\{2,3\}$-groups when $G(3)\lhd G$] \, $\\$
%Let $n=|G|$.
\indent Let $S=\{p_1^{(2)},p_2^{(2)},p_1^{(3)},p_2^{(3)}\}$ be a tame supporting set for $G$. Extend $T$ to a finite set $T_0$  disjoint from $S$ which  satisfies condition $(A_{\mQ(\mu_n)})$ where $n=|G|$.

%Assume $G$ has a normal $3$-Sylow subgroup $G(3)$ and let $S=\{p_1^{(2)},p_2^{(2)},p_1^{(3)},p_2^{(3)}\}$.
As shown in Section \ref{2-group.sec}, there exists a $\{p_1^{(2)},p_2^{(2)}\}$-adequate $G/G(3)$-extension $M/\mQ$  in which $T_0\cup\{p_1^{(3)},p_2^{(3)}\}$ splits completely.

The embedding problem $G\ra \Gal(M/\mQ)$ splits and by Lemma \ref{disjoint.lem}, $M\cap \mQ(\mu_n)=\mQ$. Thus, we may apply  \cite{Neu} and embed $M/\mQ$  into a $G$-extension $L/\mQ$ such that $T$ splits completely in $L$ and
\begin{equation*}\label{section2- normal 3-sylow condition}
\Gal(L_{p_i^{(3)}}/\mQ_{p_i^{(3)}})\cong G(3)\mbox{ for }i=1,2.
\end{equation*}
Therefore $L/\mQ$ is an $S$-adequate  $G$-extension in which $T$ splits completely, as required.
\end{proof}

%\subsubsection{Admissibility of $\{2,3\}$-groups}
%\subsubsection{When the $3$-Sylow subgroup is not normal}
%The following Proposition generalizes \cite[Theorem 1]{Son3}. Surprisingly,  the complicated proof of \cite{Son3} can be adjusted to give:
For $\{2,3\}$-groups $G$ that do not have a normal $3$-Sylow subgroup, we show that the proof of \cite[Theorem 1]{Son3} can be adjusted to give Theorem \ref{Q.thm}. %For the most part,  the proof in \cite{Son3} remains unchanged %we avoid repeating it fully here. Instead,
%\begin{prop}%\label{section2- ++splits completely proposition of 2,3}
%Let $G$ be a Sylow metacyclic $\{2,3\}$-group. Let $T$ be finite set of rational primes and $S$ a tame supporting set for $G$  disjoint from $T$. %so that $W\cap T=\emptyset$. % and in case $G(3)\lhd G$, $v_i(2)\not=3$.
%Then there exists a $G$-extension compatible with $S$ in which  $T$ splits completely.
%%(2) for every $v_i(p)\in T$, $Gal(L_{v_i(p)}/\mQ_{v_i(p)})\supseteq G(p)$ for some $p$-Sylow subgroup $G(p)$ of $G$.
%\end{prop}

%We indicate how Sonn's proof can be adjusted. Let us assume $G(3)$ is not a normal subgroup of $G$.
Let $F=F(G)$ denote the Fitting subgroup of $G$ and $F(2)$ and $F(3)$ its $2$-Sylow and $3$-Sylow subgroups, respectively. The approach of \cite{Son3} is to construct an adequate $G/F$-extension $N/\mQ$ and embed it into an adequate $G/F(2)$-extension $E/\mQ$ and an adequate $G/F(3)$-extension $L/\mQ$. Since $[EL:L]$ and $[EL:E]$ are coprime,  the compositum $EL$ is an adequate $G$-extension of $\mQ$. % $G/F(3)$ and $G/F(2)$

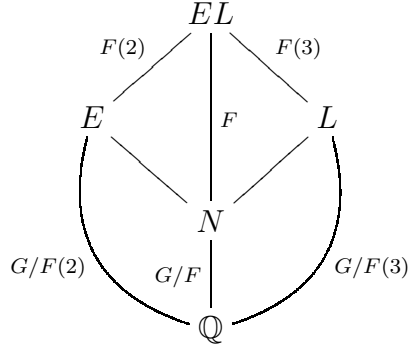
\begin{equation}\label{fields.diag}\xymatrix{
& EL \ar@{-}[dl]_{F(2)} \ar@{-}[dr]^{F(3)} \ar@{-}[dd]^{F}& \\
E \ar@{-}@/_2pc/[ddr]_{G/F(2)}\ar@{-}[dr] & & L \ar@{-}[dl] \ar@{-}@/^2pc/[ddl]^{G/F(3)}\\
&N \ar@{-}[d]_{G/F} & \\
& \mQ & }\end{equation}

When $G(3)$ is not a normal subgroup of $G$, \cite{Son3} shows that $G/F$ is isomorphic either to (1) $S_3$ or to (2) $\mZ/3$ and that in these cases we have the following partition into subcases:
$$
\begin{tabular}{| c | c | c | c | c |}
\hline
Case & $G/F$ & $F(2)$ & $G/F(3)$ & $2$-Sylow  \\
\hline
1.1 &  $\mZ/3$ & $\mZ/{2^u}\times \mZ/{2^u}$ & $\mZ/3\ltimes (\mZ/{2^u}\times \mZ/{2^u})$ &  $\mZ/{2^u}\times \mZ/{2^u}$ \\
1.2 &  $\mZ/3$ & $Q_8$ & $\SL_2(3)$ & $Q_8$\\
2.1 &  $S_3$ & $\mZ/{2}\times \mZ/{2}$ & $S_4$ & $D_8$ \\
2.2 &  $S_3$ & $Q_8$ & $S_4^*$ or $S_4^{**}$ & $Q_{16}$ or $D_{16}^*$ \\
\hline
\end{tabular}
$$
% and $F(2)$ is either the quaternion group  $Q_8$ or a group of the form $C_{2^u}\times C_{2^u}$, i.e. a homocyclic group. The following cases cover all possibilities:
%\begin{enumerate}
%\item[Case I.1.]  $G/F \cong A_3$ and $F(2) \cong C_{2^u}\times C_{2^u}$. In such a case $G/(F(3))$ is the unique
%extension of $C_{2^u}\times C_{2^u}$ by a non trivial automorphism of order $3$.
%
%\item[Case I.2.]  $G/F \cong A_3$ and $F(2) \cong Q_8$. In such a case $G/(F(3))\cong SL_2(3)$ (the unique
%extension of $Q_8$ by a non trivial automorphism of order $3$).
%%is the unique extension %of $F_2$ by a
%%non-trivial automorphism of order $3$.
%
%\item[Case II.1.]  $G/F \cong S_3$ and $F(2)$ is homocyclic. Then $F(2) \cong C_2\times C_2$ and
%$$G/(F(3)) \cong S_4.$$
%
%\item[Case II.2.] $G/F \cong S_3$, $F(2) \cong Q_8$. Then $G/(F(3))$ is one of the two
Here $Q_8$ is the quaternions group, $D_8$ the dihedral group of order $8$, $S_4^*$ and $S_4^{**}$ are the two
central extensions of $S_4$ with kernel $\mZ/2$, %, denoted by
%.
%The groups $S_4^*$ and $S_4^{**}$ have metacyclic
and %\begin{equation*}
$$\begin{array}{l}
Q_{16} = \langle x,y | x^2 = y^4, y^8=1, x^{-1}yx = y^7 \rangle, \\
%\end{equation*} \begin{equation*}
D_{16}^* = \langle x,y | x^2 = y^8 = 1, x^{-1}yx= y^3 \rangle   \\ %\end{equation*}
\end{array}$$
are their $2$-Sylow subgroups, respectively.
%\end{enumerate}

In all of the above cases the $2$-Sylow subgroups  have unique parameters $m,n$ and $t$ \footnote{The parameter $i$ is also unique up to multiplication by an odd number.}. % as described by the following lemma. % that is based on \cite{Lid2}.

\begin{lem}\label{section2- metacyclic unique representation}
Let $G \cong \mathcal{M}(m,n,i,t)$.
\begin{enumerate}
\item[(a)]  If $G\cong \mZ/{2^u}\times \mZ/{2^u}$ then $m=2^u,n=2^u,t=1$.

\item[(b)] If $G\cong Q_8$ then $m=2,n=4,t=3$.

\item[(c)] If $G\cong D_8$ then $m=2,n=4,t=3$.

\item[(d)] If $G\cong D_{16}^*$ then $m=2,n=8,t=3$.

\item[(e)] If $G\cong Q_{16}$ then $m=2,n=8,t=7$.
\end{enumerate}
\end{lem}

\begin{proof} (1) Let $x,y$ be the generators of a presentation $\mathcal{M}(m,n,i,t)$. Since $m,n|2^u$ and $mn=|G|=2^{2u}$, one has $m=n=2^u$. %Since $G$ has exponent $2^u$, one has $1=x^m=y^i$ and hence $i=0$.
For $1<t<2^u$ the group $\mathcal{M}(2^u,2^u,i,t)$ is non-abelian and hence $t=1$.

(b)--(e) are conclusions from   \cite[Theorem 22, Case $3$]{Lid2}. In this theorem, Liedahl gives necessary and sufficient conditions on a presentation $\mathcal{M}(m,n,i,t)$ of a group as in  (b)--(e) to have  an equivalent presentation with other parameters. However, these conditions require $m\geq 4$ \footnote{Note that our $m$ is denoted as $2^m$ in the notation of \cite{Lid2}.} which fails for the presentations in (b)--(e). \end{proof}

%-----------------------Proof of Sonn's theorem ---------------------

\begin{proof}[Proof of Theorem \ref{Q.thm} for $\{2,3\}$-groups when $G(3)$ is not normal] $ \\ $
\indent
We claim that the fields $N,L,E$ in diagram (\ref{fields.diag}) can be in fact chosen to be $S$-adequate extensions of $\mQ$ in which $T$ splits completely. This will imply that $EL/\mQ$ is an $S$-adequate $G$-extension in which $T$ splits completely, as required.

 We first construct the field  $E$ and let $N=E^{F/(F(2))}$.   %Let $P_3$ be a Sylow subgroup of $G$.

%Enlarge $T$ to a finite set $T_0$ which satisfies condition $(A_{\mQ(\mu_n)})$ where $n=|G|$.

In Case (1), $G/F(2)\cong G(3)$ is of odd order and therefore  \cite{Neu} gives a $\{p_1^{(3)},p_2^{(3)}\}$-adequate $G(3)$-extension $E/\mQ$  in which $T\cup\{p_1^{(2)},p_2^{(2)}\}$ splits completely. %We then let $N=E^{F(3)}$.

In Case (2), let $q\equiv 1$(mod $8$) be a prime which is not in $S\cup T$ and such that $p_1^{(2)}p_2^{(2)}q\equiv 1$ (mod $p$) for all $p\in T_0\cup \{p_1^{(3)},p_2^{(3)}\} $. Let $k=\mQ(\sqrt{p_1^{(2)}p_2^{(2)}q})$ and  let $\frak{q}$ be the prime of $k$ which lies above $q$. Note that $k$ is $\{p_1^{(2)},p_2^{(2)}\}$-adequate and $T$ splits completely in $k$. %, and %by Lemma \ref{disjoint.lem}
%the only roots of unity in $k$ are $\pm 1$.
Since the embedding problem $G/F(2)\ra \Gal(k/\mQ)$ splits, we may apply \cite{Neu} and embed $k$ into an $S$-adequate $G/F(2)$-extension $E/\mQ$ in which $T$ and $\frak{q}$  split completely.
In both Cases (1) and (2), $E/\mQ$ is $S$-adequate and $T$ splits completely in $E$.

The construction of the field $L$ is the same as in \cite{Son3} with few modifications. Since the construction in \cite{Son3} is involved and long, we do not repeat it here. A self contained version of the modified construction  can be found in the author's thesis (\cite{Nef}). For the reader to whom \cite{Son3} is available,  we give below the list of required modifications.

% that are based on using Corollary \ref{neu.gen}.
Note that our field $N$ was denoted in Case (1) of \cite{Son3}  by $k$ and in Case (2) by $K$. %Since our $N$ differs from the construction in \cite{Son3}, we indicate how the construction of the field $L$ in \cite{Son3} should be adjusted: %which is compatible with $S$ and in which $T$ splits completely.
%Since case (2.1) follows from case (2.2), we remain with three cases.
%With the above definitions of $K$ and $E$, Sonn's construction of the field $L$ applies with the following changes.
\begin{enumerate}
\item Replace the primes $p_1,p_2$ (resp. $p,q$) in Case (1) (resp. Case (2)) of \cite{Son3}  by the primes $p_1^{(2)},p_2^{(2)}$, respectively. Since $S$ is a supporting set, Lemma \ref{section2- metacyclic unique representation} implies that the prime $p_1^{(2)},p_2^{(2)}$ satisfy the congruence relations required in \cite{Son3}  from  $p_1,p_2,p,q$. Note that in Case (1), the primes $p_1^{(2)},p_2^{(2)}$ split completely in $N(\mu_{n})$ as required in \cite{Son3}. %Also note that we assumed that primes in a tame supporting set are greater than $3$ since this assumption is required from $p,q$ for the proof of case 2 in \cite{Son3}.%
    Also note that since  $p_1^{(2)},p_2^{(2)}$ are prime to $|G|$, they are greater than $3$ as required in Case (2) of \cite{Son3}.
\item In Case (2), we add the prime $\frak{q}$ to the modulus $\frak{m}$ and require that the element $\gamma$ is congruent to $1$ mod $\frak{q}$. The field $M$ constructed in \cite{Son3} then satisfies  $\Gal(M_q/\mQ_q)\cong \mZ/2\mZ$.
 Since $q\equiv 1$ (mod $8$), $\Gal(M_q/\mQ_q)$ can be embedded into a $\mZ/4\mZ$ extension and therefore
 %and since $\frak{q}$ splits completely in $M$,
 the embedding problem $G/(F(3))\ra \Gal(M/\mQ)$ is solvable at $q$ as well. As shown in \cite{Son3} it is solvable at all other primes and hence globally solvable. %Indeed, the kernel of the embedding problem is $\mZ/2\mZ$ and $\Gal(M_q/\mQ_q)\cong \mZ/2\mZ$ embeds into a $\mZ/4\mZ$ extension. Therefore one can proceed with the same argument to show that $G/(F(3))\ra \Gal(M/\mQ)$ has a global solution.
\end{enumerate}
With these changes the field $L$ constructed in \cite{Son3} gives an $S$-adequate $G/F(3)$-extension. In order to make the set $T$ split completely in $L$, we make the following additional modifications:
\begin{enumerate}
\item In Case (1.1),  the embedding problem $G/(F(3))\ra \Gal(N/\mQ)$ splits and hence has the trivial solution $\phi$. Instead of applying Theorem \ref{7cond.thm},
we apply Corollary \ref{Neu.cor} insuring the same prescribed conditions at $S$ but in addition that the solution remains trivial at primes of $T$.
\item In Cases (1.2) and (2), we add the primes of $N$ that lie over primes of $T$ to the modulus $\frak{m}$ and require that $\gamma\equiv 1$ (mod $\fp$) for every $\fp\cap\mQ\in T$. This insures that $T$ splits completely in $K$ (resp. in $M$) in Case (1.2) (resp. Case (2)).
%\item Denote the field $K$ (resp. $M$) appearing in \cite[Case 1.2]{Son3} (resp. \cite[Case 2]{Son3}) by $M$. Then (2) insures that  $T$ splits completely in $M$.
%\item

Let $\phi$ be the solution obtained in Case (1.2) (resp. Case (2)) of \cite{Son3} to the embedding problem $G/F(3)\ra \Gal(K/\mQ)$ (resp. $G/F(3)\ra \Gal(M/\mQ)$).
We apply Theorem \ref{7cond.thm} in order to change $\phi$ to a solution $\psi$ which is trivial at primes of $T$. Since the local embedding problem at $p_i^{(2)}$ is Frattini, $\psi$ is surjective at $p_i^{(2)}$, $i=1,2$. Thus, the fixed field $L$ of $\ker\psi$ is $S$-adequate and $T$ splits completely in $L$.
\end{enumerate}
\end{proof}
\newpage
\bibliographystyle{plain}

\end{document}